\documentclass{elsart1p}

\usepackage{amsmath,,amssymb,amsfonts}
\usepackage[american]{babel}
\usepackage{cite}
\allowdisplaybreaks

 \newtheorem{theorem}{Theorem}[section] 
 \newtheorem{lemma}[theorem]{Lemma} \theoremstyle{definition}
 \newtheorem{definition}[theorem]{Definition}
 \newtheorem{proposition}[theorem]{Proposition}
 \newtheorem{corollary}[theorem]{Corollary}
 \newtheorem{obs}[theorem]{Remark}

\newcommand{\vp}{\vspace{0,5cm}}   
\newcommand{\vi}{\vspace{0,2cm}}   
\newcommand{\N}{{\mathbb{N}}}   
\newcommand{\Z}{{\mathbb{Z}}}   
\newcommand{\Q}{{\mathbb{Q}}}   
\newcommand{\U}{{\mathcal{U}}}   
\newcommand{\M}{{\mathcal{M}^{0}}}   
\newcommand{\m}{{\mathcal{M}}}   
\newcommand{\A}{{\mathcal{A}}}  
\newcommand{\I}{{\mathfrak{I}}}  
\newcommand{\z}{{\mathcal{Z}}}  
\newcommand{\s}{{\mathcal{S}}}  
\newcommand{\E}{{\mathcal{E}}}  
\newcommand{\dg}{{\mathbf{diag}}}  

\begin{document}

\begin{frontmatter}


 \title{Hyperbolicity of Semigroup Algebras}

\author[UFABC]{E. Iwaki},
\ead{edson.iwaki@ufabc.edu.br}
\author[IME]{S. O. Juriaans},
\ead{ ostanley@ime.usp.br}
\author[IME]{A. C. Souza Filho\corauthref{cor}\thanksref{now}}
\corauth[cor]{Corresponding author.}
\ead{calixto@ime.usp.br}
 \thanks[now]{This article corresponds to the second chapter of the third author's
PhD. Thesis, see \cite{thss}}
\address[UFABC]{Centro de Matem\'atica, Computa\c c\~ao e Cogni\c c\~ao, Universidade
 Federal do ABC, Rua Catequese, 242, $3^{\circ}$ andar, 09090-400, Santo
 Andr\'e, SP, Brasil}
\address[IME]{Instituto de Matem\'atica e Estat\'\i stica,
 Universidade de S\~ao Paulo, Caixa Postal 66281, S\~ao Paulo, CEP
 05315-970 - Brasil}


\title{}


\author{}

\address{\scriptsize Mathematics subject Classification Primary [$16U60$, $20M25$]; Secondary [$20M10$, $16S36$]}

\begin{abstract}
Let  $\A$ be a  finite dimensional
$\Q$-algebra and $\Gamma \subset \A$ a $\Z$-order. We classify those $\A$ with
the property that ${\Z^{2} \not
\hookrightarrow} \U (\Gamma)$. We call this last property {\it the hyperbolic
property}. We apply this in the case that $\A=K S$ a semigroup algebra with
$K=\Q$ or $K=\Q(\sqrt{-d})$. In
particular, when $K S$ is semi-simple
  and has no nilpotent elements, we prove that
  $S$ is an inverse semigroup which is the disjoint union of Higman
  groups and at most one cyclic group $C_n$ with $n \in \{5,8,12\}$.
  \end{abstract}

\begin{keyword} {Semigroup, Semigroup Algebras, Hyperbolic Groups, Group Rings, Units }


\end{keyword}

\end{frontmatter}

\section{Introduction}

\vp

In this paper we focus on what we call the {\it hyperbolic
property}. We say that a finite dimensional $\Q$-algebra  $\A$
has the hyperbolic property if for every $\Z$-order $\Gamma \subset
\A$ the unit group $\U(\Gamma)$ does not contain a finitely
generated abelian group of rank greater than one. This terminology
is suggested by the fact that hyperbolic groups have this property
\cite{grmv}.

\vi

Research in this direction goes back to Jespers, who classified
those finite groups $G$ for which $\U(\Z G)$ has a non-Abelian free
normal complement \cite{jsp}. More recently Juriaans-Passi-Prasad
have given contributions on this topic in the integral group ring
case \cite{jpp}, and Juriaans-Passi-Souza Filho in the group ring
$RG$ when $R$ is the ring of algebraic integers of a quadratic
rational extension \cite{jpsf}.

\vi

Here we give a complete classification of
the finite semigroups whose semigroup algebra $KS$ has the hyperbolic property with $K=\Q$ or $k=\Q(\sqrt{-d})$. Part of this
was done by  Jespers and Wang \cite{gdf} who classified the finite
semigroups $S$  for which the unit group $\U(\Z S)$ of the integral
semigroup ring $\Z S$ (we of course assume that this ring contains
an identity) is a finite group.
Firstly, we prove a structure theorem for the finite dimensional
$\Q$-algebras with the hyperbolic property. We prove that the
radical of such an algebra has nilpotency index at most $2$ and that its Wedderburn-Mal\^cev
components consist of copies of $\Q$ or quadratic fields, totally
definite quaternion algebras, two-by-two matrices over $\Q$
and upper-triangular matrices over $\Q$. Details on the
structure of these algebras are given in section $3$. In  section
$4$, we classify the finite semigroups $S$ whose semigroup algebras
$K S$  has the hyperbolic property, with $K=\Q$ or $K$ is a
quadratic extension of $\Q$. In section $5$ we study the idempotents
of the maximal subgroups of finite semigroups $S$ which are not
semi-simple in order to obtain a best comprehension of the structure
of $S$, when $\Q S$ has the hyperbolic property.

\vi

Notation is mostly standard and we refer the reader to  \cite{clfrd}
and \cite{jok} for the theory of semigroup and semigroup algebra.
However, for the readers convenience, section $2$ contains
some basic facts on the theory of semigroups.

\vp

\section{Preliminaries}

\vp

 A non-empty set $S$  with an associative binary
operation $\cdot:S^{2} \rightarrow S$ is a semigroup. Let $S$ be a
semigroup, the set $S^{1} \dot{=} S \cup \{1\}$, such that, $\forall
s \in S, s \cdot 1=1 \cdot s=s$ is a monoid, that is, a semigroup with an
identity element, and the set $S^{\theta} \dot{=}S \cup \{\theta\}$,
such that, $\forall s \in S, s
\cdot \theta=\theta \cdot s=\theta$, $\theta$ called a zero element, is a semigroup with
a zero element. A semigroup $S$ with zero $\theta$ is a null semigroup
if for all $x,y \in S$, $x \cdot y=\theta$. An element $e \in S$, such that,
$e^{2}=e$ is an idempotent. Denote $E(S):=\{e \in S/e^2=e\}$
 the set of idempotents
of $S$ and let $e,f \in E(S)$; then $e \leq f$ if $e \cdot f=f\cdot
e=e$. An idempotent $f \in E(S)$ is primitive if $f \neq \theta$ and
if $e \leq f$ yields $e=\theta$ or $e=f$. A semigroup $S$ is simple,
if it does not properly contain any two-sided ideal. A semigroup
with zero $\theta$ is $0$-simple if $S^2 \neq \{\theta\}$ and
$\{\theta\}$ is the only proper two-sided ideal of $S$. A $0$-simple semigroup $S$ is
completely $0$-simple if it contains primitive idempotents.

\vi

Let $I$ be an ideal of $S$. The semigroup of the Rees factors,
denoted by $S/I$, is the set $(S \setminus I) \cup \{\theta\}$
subject to the operation $\cdot$ defined by $$s \cdot t := \left \{
\begin{array}{ll} s \cdot t,& \textrm{ if } s \cdot t \notin I \\ \theta,&
\textrm{ if } s \cdot t \in I \end{array} \right.$$

\vi

 A principal series of a
semigroup $S$ is a chain
$$S=S_{1}\supset S_{2}\supset \cdots \supset S_{n}\supset
S_{n+1}=\emptyset$$ of ideals $S_{i}$ of $S$, such that $S_{i+1}$ is
a maximal ideal of $S_{i}$, $1 \leq i \leq n$. The semigroups of the
Rees factors $S_{i}/ S_{i+1}$ are called factors of this principal
series. It is well known that if $S$ is a finite semigroup, then the
factors of $S$ are either a null semigroup with two elements which
we will call null factor or a completely $0$-simple semigroup.

\vi

The semigroups $S$ which are union of groups appear naturally in the
context we work. Since a semigroup which is a union of groups  is the
disjoint union of its maximal subgroups we have the following:

\vi

\begin{lemma}
 \label{GD}
 Let $S$ be a finite semigroup whose factors are isomorphic to groups with a
 zero element adjoined
$\theta_i$, that is, $S_{i}/ S_{i+1} \cong G_{i}^{\theta_i}$. Then
$S$ is a disjoint union of groups.
\end{lemma}

\vi

Let $G$ be a group and $I$ and  $\Lambda$ arbitrary non-empty sets.
By an  $I \times
 \Lambda$ Rees matrix, we mean an $I \times \Lambda$ matrix over $G^{\theta}$
 with at most a unique entry different from $\theta$. For $a \in G, i \in I$ and $\lambda \in
 \Lambda$, $(a)_{i \lambda}$ denotes an $I \times
 \Lambda$ Rees matrix over $G^{\theta}$, where $a$ is the entry corresponding   to   row $i$
 and   column $\lambda$ and all other entries are zero. For any $i \in I$ and $\lambda \in
 \Lambda$, the expression $(\theta)_{i \lambda}$ denotes the $I \times
 \Lambda$ null matrix, which is  also denoted  by $\theta$.

\vi

Since we are dealing with finite semigroups, it is sufficient to
consider a finite number of rows and columns, $m,n$, respectively.
  For $1\leq i \leq m$ and $1 \leq j\leq n$, fix $P=(p)_{ij}$ a $m \times n$
 matrix over $G^{\theta}$, called a sandwich matrix, and let $\M$ be the set
 of the  $m \times n$ Rees matrices over  $G^{\theta}$.
 In $\M$ we define the operation
$AB =A \circ P \circ B,$ where $\circ$ denotes the usual matrix
product, which is binary and associative and therefore the set
$\{\M,\circ\}$ is a semigroup. This semigroup is denoted by
 $\M(G;m,n; P)$ and $G$ is called its structural group.

\vi

In a similar way we define the Munn  matrices. Let $R$ be a ring and
$m,n$ positive integers. Consider $\m(R;m,n; P)$ the set of $m
\times n$ matrices over $R$. For each $A=(a_{ij}),B=(b_{ij})
 \in \m(R,m,n, P), 1\leq i \leq m, 1 \leq j\leq n$,   addition is defined by $A+B=(a_{ij}+b_{ij})$, and
multiplication by $AB=A \circ P \circ B$, where $P$ is a fixed $n
\times
 m$ matrix with  entries in $R$ and $\circ$ is the usual matrix operation.
 The ring $\m(R;m,n; P)$ is called an algebra of matrix type over $R$ or a matrix algebra over $R$.

\vi

Let $\A$ be a finite dimensional $\Q$-algebra. A
 unitary subring $\Gamma$ of $\A$ is called a $\Z$-order, or simply
 an order, if it is a finitely generated $\Z$-submodule such that
 $\Q \Gamma=\A$, (see \cite[$1.4$]{uni}). Remember that, by the
Borel-Chandra Theorem \cite{bc}, the unit group of a $\Z$-order of
$\A$
 is finitely generated and hence the hyperbolicity of $\U(\Gamma)$ makes sense.
Furthermore, if $\U(\Gamma)$ is a hyperbolic group then
$\U(\Gamma_0)$ is hyperbolic for all $\Z$-order $\Gamma_0 \subset
\A$, since the unit groups of orders are commensurable. It is known,
\cite{grmv}, that the
 hyperbolicity of $\U(\Gamma)$ implies that  ${\Z^{2} \not
\hookrightarrow} \U (\Gamma)$. This suggests the following
definition.

\vi

\begin{definition} Let $\A$ be a finite dimensional $\Q$-algebra and
 $\Gamma$ a $\Z$-order of $\A$. We say that $\A$ has the
hyperbolic property if $\Z^{2} {\not \hookrightarrow} \U (\Gamma).$
\end{definition}

\vi

Note that, as seen above, this definition does not depend on the
particular order $\Gamma$ of  $\A$, (see \cite{bc}).

\vi

Throughout the text we use the standard notation $\dg(a_1, \cdots,
a_n)$ for a $n \times n$ matrix with elements on the main diagonal
set to $a_1, \cdots, a_n$ and all the other elements set to zero.
Also $e_{ij}$ denotes the elementary matrix whose entry is $1$ in
the position $i,j$ and zero otherwise. We denote by  $T_2 (\Q
):=\left ( \begin{array}{ll}\Q & \Q \\0 &\Q \end{array} \right )$
with the usual matrix multiplication.

\vp


\vp

\section{Finite Dimensional Algebras with the Hyperbolic Property}

\vp

The main result of this section is Theorem \ref{tfadfh} in which we
classify the finite dimensional $\Q$-algebras which have the
hyperbolic property.

\vp
\begin{theorem}\label{tfadfh}Let $\A$  be a finite dimensional $\Q$-algebra,
$\A_i$ a Wedderburn component of $\A$ and $\Gamma_i \subset \A_i$ a $\Z$-order.
 Then
\begin{description}
\item[$(i)$]  $\A$ has the hyperbolic property, is semi-simple and without
  nilpotent elements if, and only if,
  $$\A =\oplus \A_{i},$$ and
  for at most one index $i_{0}$ we have that  $\U(\Gamma_{i_{0}})$ is hyperbolic and
  infinite.
\item[$(ii)$]  $\A$ has the hyperbolic property and is semi-simple with
  nilpotent elements if, and only if,
$$\A =(\oplus \A_{i})\oplus M_2(\Q).$$
\item[$(iii)$]  $\A$ has the hyperbolic property and is non-semi-simple with
  central radical if, and only if, $$\A =(\oplus \A_{i})\oplus J.$$
\item[$(iv)$]   $\A$ has the hyperbolic property and is non-semi-simple with
  non central radical if, and only if, $$\A =(\oplus \A_{i})\oplus
 T_{2}(\Q).$$
\end{description}
For each item $(i)-(iv)$, the $\A_{i}$'s are either at most  a
 quadratic imaginary
 extension of $\Q$ or a totally definite quaternion algebra.
Furthermore, in the decompositions in $(i)-(iv)$ above the direct
summands are ideals.
\end{theorem}

\vi

We will consider $\A$ a finite dimensional  $\Q$-algebra with
radical $J(\A)$. According to a theorem of Wedderburn-Mal\^cev,
there exists a semi-simple subalgebra $\s (\A)$ of $\A$ such that
$$\A = \s(\A) \oplus J(\A), \textrm{ as a  direct sum of vector spaces},$$
with $J(\A)$ the Jacobson radical.

\vi

\begin{lemma} \label{rsmsp} Let ${\A=\s(\A) \oplus J(\A)}$ be a finite dimensional $\Q$-algebra
with the hyperbolic property. Then $J^{2}=0$, $dim_{\Q}(J(\A))\leq
1$, that is, as a vector subspace  of $\A$, $J=J(\A)$ has  dimension
at most $1$. Furthermore, if $J \neq 0$, there exists $j_0 \in \A$
such that $j_0^2=0$, ${J=\langle j_0\rangle_{\Q}}$, and $1+J \cong
\Q$ as a multiplicative group.
\end{lemma}

\begin{pf}
Obviously $1+J$ is a multiplicative torsion free  nilpotent group.
Let $G$ be any finitely generated subgroup of $1+J$. Hence $\z(G)
\neq 1$. Since $\Z^2 {\not \hookrightarrow} (1+J)$ the same holds
for $G$, so $G=\z(G) \cong \Z$. Since $J$ is a  nilpotent ideal,
there exists a least positive integer $n$, $J^n=0=JJ^{n-1}$ thus
$1+J=\z(1+J)=\z(1+J^{n-1})\cong \Z$. Hence $n=2$ and $dim_{\Q}J \leq
1$. If $J \neq 0$ and $\Gamma \subset \A$ is a $\Z$-order, let $x,y
\in J \cap \Gamma$. Then the group $\langle 1+x,1+y\rangle <
\U(\Gamma)$, and $1+x,1+y$ are units of infinite order. Since
$\U(\Gamma)$ is hyperbolic we have $\langle 1+x,1+y\rangle \cong \Z
$. Hence $\langle 1+x \rangle \cap \langle 1+y \rangle$ is
non-trivial and there exist integers $ m,n$, such that,
$(1+x)^{m}=(1+y)^{n}$. Since $x,y$ are $2$-nilpotent, we have
$1+mx=1+ny$, and thus $x=\frac{n}{m}y$. So $\{x,y\}$ is a
$\Q$-linearly dependent set and we conclude that $dim_{\Q}(J)=1$.
Write $J=\Q j_0$, so $1+J \cong \Q$. Indeed, $\phi:1+J \rightarrow
\Q, \phi(1+qj_0)=:q$ is an isomorphism. $\square$
\end{pf}

If $\s_1(\A)$ and $\s_2(\A)$ are subalgebras of $\A$ such that $\A =
\s_l(\A) \oplus J(\A), l=1,2$, then there exists $n \in J(\A)$ such
that $S_1=(1-n)\s_2(\A)(1-n)^{-1}$, hence the semi-simple subalgebra
$\s(\A)$ is unique, up to isomorphism \cite[Theorem $72.19$]{cr}.
Choose a semi-simple subalgebra $\s(\A)$ of $\A$ and let $\E(\A)=\{\
E_{1}, \cdots, E_{N}\}, N \in \Z^{+}$ be the set of orthogonal
central primitive idempotents of the semi-simple subalgebra $\s(\A)$
such that \linebreak $\A_i:=\s(A)E_i$, we have $\s(\A)\cong \oplus
\A_i$. The algebra $\A$ can be written as a finite direct sum of
indecomposable two-sided ideals, \cite[Theorem $55.2$]{cr}. Let
$$\A=(\bigoplus_{E_{i} \in \E(\A)} \s(\A) E_{i}) \oplus J(\A)$$
be its Wedderburn-Mal\^cev decomposition.

\vp

\begin{proposition}\label{pcent}
Let $\A$ be a finite dimensional non-semi-simple $\Q$-algebra with
\linebreak $dim_{\Q}(J(\A))=1$, $J(\A)=\langle j_0 \rangle$ and
$N=|\E(\A)|$. The following conditions hold:
\begin{enumerate}
\item For all $x \in \A$, there exist
 $\lambda_{x},\mu_{x} \in \Q$, such that, $xj_0=\lambda_{x}j_0$ and $j_0x=\mu_{x}j_0$.
\item If $x$ is an idempotent, then $\lambda_x,\mu_{x} \in\{0,1\}$.
\item There exist unique $E,F \in \E(\A)$, such that $Ej_0\neq0$ and $j_0F\neq0$.
\item If $E=F$ then $J$  is central.
\item If $J$ is non-central then, up to an index reordering, we can suppose
that
 $E=E_{1}$, and $F=E_{N}$ and  $E_1j_0=j_0E_N=j_0$. In particular, in this
 case $N \geq 2$.
\end{enumerate}
\end{proposition}

\begin{pf}
Since $J$ is a two-sided ideal of $\A$, for $x \in \A$
we have  $x \cdot j_0 \in \langle j_0 \rangle_\Q$, hence there
exists  $\lambda_x \in \Q$, such that,  $x \cdot j_0=\lambda_x j_0$.
Similarly, for some $\mu_x \in \Q$, $j_0 \cdot x=\mu_x j_0$.

If $x$ is  an idempotent then $x \cdot j_0=\lambda_x j_0 =x^2 \cdot j_0=x\cdot(x
\cdot j_0)=x\cdot (\lambda_x j_0)= \lambda_x^2 j_0 \cdot x$ and hence
$(\lambda_x^2-\lambda_x)j_0=0$. It follows that $\lambda_x^2-\lambda_x=0$ and
so $\lambda_x \in \{0,1\}$, in a similar way $\mu_x \in \{0,1\}$.

Since $ 1=\sum_{1\leq i \leq N}E_{i}$ we have that  $ 1
  \cdot j_0=\sum_{1\leq i \leq
 N}(E_{i} \cdot j_0)=\sum_{1\leq i \leq N}(\lambda_ij_0)=(\sum_{1\leq i \leq N}\lambda_i)j_0$,
 and so, $  \sum_{1\leq i \leq N}\lambda_i=1$. Since each $E_i$
  is an idempotent, we have that  $\lambda_i \in \{0,1\}$. This
  yields the existence of a unique index,   $m$ say, $1 \leq m \leq N$, such
  that,  $\lambda_{m}=1$ and $\lambda_i=0$ if $i\neq m$. Similarly
  there exists a unique $ k, 1 \leq k \leq N$, such that, $j_0 \cdot E_{k}=j_0$ and
    $j_0 \cdot E_i=0$ for $i \neq k$.

 If $E=F$, $E_m \cdot j_0=E_k \cdot j_0$ then, by  uniqueness, $m=k$ and
 $j_0 \cdot E_m=E_m \cdot j_0=j_0$. On the other hand, for  $i \neq k$, we have that
  $j_0 \cdot E_i =E_i \cdot j_0 =0$. Therefore $J$ commutes with $\s(\A)$ and thus it is
 central. The other statements are now also clear.
$\square$ \end{pf}

\begin{corollary} \label{jcen}
Let $\A$ be a finite dimensional non-semi-simple $\Q$-algebra
with the hyperbolic property. Then
 $J(\A)=\langle j_0\rangle$ is central in $\A$ if, and only if, there exists a
 unique $E \in \E(\A)$ such that $j_0=E \cdot j_0=j_0 \cdot E \neq 0$.
\end{corollary}

\vp

\begin{obs}\label{nhl} Denote by $M$ the left annihilator of $J$ in $\A_{1}$. Since
$dim_{\Q}(J) =1$, it follows that $M$ is a proper ideal of $\A_{1}$.
 If $x \in \A_{1}$, then there exists $\lambda_{x} \in \Q$, such that $x \cdot
j=\lambda_{x}j_0$, so $x=(x-\lambda_{x}E_1)+\lambda_{x}E_1$ and thus
$\A_{1}=M \oplus \Q E_{1}$ which implies that $dim_{\Q}(M)+1 =
dim_{\Q}(\A_{1})$. Since $\A$ is a simple $\Q$-algebra it must
be that $M=\{0\}$ and $dim_{\Q}(\A_1) =1$.
 Analogously, we obtain $\A_{N} \cong \Q$.
\end{obs}

\vp

 Let $\A$ be a rational finite dimensional algebra with the hyperbolic
 property, $N=|\E(\A)|\geq 2$, and
${E_{1},E_{N} \in \E(\A)}$ the idempotents with the property
$E_{1} \cdot j_0=j_0 \cdot E_{N}=j_0$. We have:
$$\s(\A)= \underbrace{(\oplus_{1< i < N} \s(\A) \cdot E_{i})}_{B} \oplus
\underbrace{(\s(\A) \cdot E_{1} \oplus \s(\A) \cdot E_{N})}_{C}=B \oplus C.$$
Letting   $\A_{i} = \s(\A) \cdot E_{i},1 \leq i \leq N$, we have:
$\A \cong B \oplus \A_{1} \oplus \A_{N} \oplus J(\A)$.
Considering the map  $$\begin{array}{llll}
\varphi:&\A_{1} \oplus \A_{N} \oplus J(\A) &\rightarrow &
T_2(\Q)
\\& a_{1} E_{1}+a_{N} E_{N}+qj_0& \mapsto & \dg(a_1,a_N)+q e_{12} \end{array}.$$
Clearly $\varphi$ is an algebra isomorphism. Hence $\A_{1}\oplus \A_{N} \oplus J(\A) \cong T_2(\Q)  .$

\vi
Thus we proved the next theorem:
\vi

\begin{theorem} \label{triangsup}
Let  $\A$ be a finite dimensional non-semi-simple $\Q$-algebra
with the hyperbolic property. If
 $J(\A)=\langle j_0\rangle$  is non-central, then, up to a reordering, we have
 $E_{1} \cdot j_0=j_0 \cdot E_{N}=j_0$, and $E_{N} \cdot j_0=j_0 \cdot E_{1}=0$. For the others idempotents
 $E_{i}, i
\notin \{1,N\}$ we have that $ E_{i} \cdot j_0=j_0 \cdot E_{i}=0$. Moreover,
$$\A_{1} \oplus \A_{N} \oplus J(\A)\cong T_2(\Q)$$ is an ideal of $\A$.
\end{theorem}

\vp

\begin{corollary} \label{sid}
Let $\A$ be a finite dimensional $\Q$-algebra with the
hyperbolic property, and $|\E(\A)|=N$. $J$ is non-central if and only
if
$$\A \cong B \oplus T_{2}(\Q) \cong \dg(\A_2, \cdots, \A_{N-1}) \oplus
T_2(\Q)  \cong  \dg(\Q,\A_2, \cdots, \A_{N-1},\Q)\oplus \Q e_{1N}.$$

Moreover, for each $1\leq i\leq N$, $\A_i$ is at most a quadratic
imaginary
 extension of $\Q$, or
 a totally definite quaternion algebra.
\end{corollary}

\vp

\begin{pf}
By the previous theorem, $B$ and  $T_{2}(\Q)$ are ideals whose
direct sum
 equals $\A$. Consider the algebra isomorphism
$$\begin{array}{llll}\varphi:&\dg(\A_2, \cdots, \A_{N-1})\oplus T_2(\Q)
&\longrightarrow & \dg(\Q, \A_2, \cdots, \A_{N-1},\Q)+\Q e_{1N}\\
&\dg(a_2, \cdots, a_{N-1})\oplus \dg(q_1,q_N)+q e_{12}
&\mapsto & \dg(q_1, a_2, \cdots, a_{N-1},q_n)+q e_{1N}
 \end{array}.
$$

Let $\Gamma_{0}  \subset \A$ be a $\Z$-order, with $E_{1}+E_{2}+
\cdots E_{N}+j_0 =1+j_0 \in \U(\Gamma_{0})$. By Lemma \ref{rsmsp},
 $\langle 1+j_0 \rangle \cong \Z$.
Suppose $\gamma_{i_0}E_{i_0} \in \Gamma_{i_0}$ is an   element of
infinite order and set \linebreak $\gamma=E_{1}+
\cdots+\gamma_{i_0}E_{i_0}+ \dots+ E_{N}$, $1 \neq i_0 \neq N$. We
have $o(\gamma ) =\infty$, $\langle 1+j_0 \rangle \cong \Z$ and
$\langle 1+j_0 \rangle \cap \langle \gamma \rangle =\{1\}$. Since
$\A_{i_0} \subset C_{\A}(J)$, the centralizer of $J \subset\A$, we
have that $\langle 1+j_0 \rangle \times \langle \gamma \rangle \cong
\Z^{2}$ is a subgroup of $\U(\Gamma_{0})$, a contradiction.
Therefore $\U(\Gamma_{i})$ is a torsion group and hence is finite.
Obviously $|\U(\Gamma_{1})\cong \U(\Gamma_{N})| \leq 2$  since, by
the previous theorem,  $\A_{1} \cong \A_{N} \cong \Q$. Therefore, by
Lemma $2.3$ of \cite{uni}, each $\A_i$ is at most an
  imaginary extension of $\Q$, or
 a totally definite quaternion algebra. To prove the converse it is enough to consider the right and left action of $J$
 on the semi-simple part of $\A$.
$\square$ \end{pf}

If $\A$ has the hyperbolic property and the radical  $J\neq \{0\}$ is central,
then   $\s(\A)$  is a direct  sum of division rings: in fact,
 if any component of $\s(\A)$ were of matrix
 type it would have an element of infinite order.
  Hence once again we could embed ${\Z^{2} \hookrightarrow \U(\Gamma)}$,
  for some $\Z$-order $\Gamma$. Therefore  the simple components  $\A_{i}, 1
\leq i \leq N$, of $\s(\A)$ are as in the corollary above.

\vp

\begin{corollary}\label{ohp}
Let $\A$ be a  finite dimensional $\Q$-algebra  with the
hyperbolic property, and $|\E(\A)|=N \geq 2$. If $J$ is central then
$\U(\Gamma_{i})$ are  finite subgroups. Moreover, if  $\Gamma
\subset \A$ is a  $\Z$-order then $\U(\Gamma)$ is commensurable with
$\Z \times C_{2} \times C_{2} \times \prod H_i$, $|H_i|<\infty$. In
particular, $\U(\Gamma)=\Phi(\U(\Gamma))$, the finite conjugacy
center of $\U(\Gamma)$.
\end{corollary}

\vi

\begin{obs}\label{hphi}
Let $A$ be a ring of characteristic zero and let
$\theta_{1},\theta_{2} \in A$ be commuting $2$-nilpotent elements. If
$\{\theta_{1},\theta_{2}\}$ is $\Z$-L.I. then $\Z^{2}$ embeds into $\U (A)$.
\end{obs}

\vi

Now we are ready to proof the main result of this section.

\vi

\begin{pf}(of Theorem \ref{tfadfh})
 Items  $(iii)$ and $(iv)$
follow from Theorem \ref{triangsup} and its corollary.

 We now prove  $(2)$: as $\A$ is semi-simple with nilpotent elements  we have
 that \linebreak  $\A \cong
\oplus M_{n_{i}}(D_{i})$, where the  $D_{i}$'s are division rings.
Remark \ref{hphi} implies that   $n_{i} \leq 2, \forall i$. The
hyperbolicity hypothesis implies that there is at most one component
with $n_{i_0}=2$ and  it is isomorphic to $M_{2}(\Q)$ (this follows
by Remark \ref{hphi}). Let $\Gamma_i \subset \A_i$ be a $\Z$-order
of $\Gamma_i$ and consider the $\Z$-order  $\Gamma_{0} = M_{2}(\Z)
\oplus (\oplus_{i\neq i_0} \Gamma_{i}) \subset \A$. We have that
$\U(\Gamma_{0}) \cong GL_{2}(\Z) \times (\prod \U(\Gamma_{i}))$. It
follows that all
 $\U(\Gamma_{i})$ are torsion  groups  and hence they are finite.

The converse is straightforward, since $GL_2(\Z)$ is  hyperbolic.

We now prove $(i)$:
 if $\A$ is semi-simple with no nilpotent elements
then $M_{2}(\Q)$ is not a Wedderburn component of $\A$ and hence $\A
\cong \oplus \A_i$, a direct sum of division rings. If for any
$\Z$-order $\Gamma \subset \A$ it holds that $\U(\Gamma)$ is finite
we are done. Suppose $|\U(\Gamma)|=\infty$, Let $\Gamma = \oplus
\Gamma_i$; then $\U(\Gamma) \cong \prod \U(\Gamma_i)$. The
hyperbolicity of $\U (\Gamma )$ implies that there can be at most
one index $i_0$ for which $\U (\Gamma_{i_0})$ is infinite and hence we are
done.   The converse is obvious. $\square$ \end{pf}

\begin{obs}\label{more}
Let $\A$ be a finite dimensional non-semi-simple  $\Q$-algebra  with the
hyperbolic property and $J$ its radical. If $a \in \A$ is a
non-trivial nilpotent element then $a \in J$.

 In fact, by Theorem \ref{tfadfh}, $\A \cong B \oplus T_{2}(\Q)$ (respectively $\A
 \cong B \oplus J$) if $J$ is non-central, (respectively if $J$ is central).
 Since each $\A_{i}, 1<i<N,$
 is a division ring, hence $a \in T_2 (\Q)$, (respectively $a \in J$).
 It is sufficient to consider the case for $J$  non-central. Let
 $a=\dg(x,z)+y e_{12}$; $a^{2}=0$ yields $x=z=0$, and $y \in \Q$. Therefore,
 $a=y e_{12} \in J$.
\end{obs}

\vi

\section{Semigroup Algebras with the Hyperbolic Property}

 \vi
 In this section we classify the finite semigroups $S$ for which $\Q
S$ has the hyperbolic property, we also classify the extensions $K=\Q
 (\sqrt{-d})$ with this property. First some terminology: a finite
group $G$ is called a {\it Higman group} if $G$ is   either abelian
of exponent dividing $4$ or $6$ or a Hamiltonian $2$-group. Recall
that nilpotent free means the absence of nilpotent elements and
$Q_{12} \cong C_3 \rtimes C_4$ where $C_4$ acts by inversion on
$C_3$.

\vi

 Let
$K$ be a field and $S$ a semigroup.  By the semigroup algebra $KS$
of $S$ over $K$ we mean an algebra $\A$ over $K$ which contains a
subset $\overline{S}$ that is a $K$-basis and a multiplicative
semigroup of $\A$ isomorphic to $S$. Let $S$ be a semigroup with a
zero element $\theta$. By the contracted semigroup algebra $K_{0}S$
of $S$ we mean an algebra over $K$ with  a basis $B$, such that, $B
\cup \{\theta\}$ is a subsemigroup of $K_{0}S$ isomorphic to $S$. If
$S$ is a Rees  matrix semigroup, $S= \M(G;m,n;P)$, then the
contracted algebra $K_{0}S \cong  \m(KG;m,n;P)$, \cite[Lemma
$5.17$]{clfrd}.

\vi

We suppose the algebra $KS$ has a unity. By \cite[Corollary
$5.26$]{jok}, if $S=\M(G;m,n;P)$ is a Rees matrix semigroup then the
following conditions are equivalent:
\begin{enumerate}
 \item   $\Q_{0}S$ is unitary;
 \item $m=n$ and $P$ is an invertible matrix in $M_{n}(\Q G)$.
\end{enumerate}

\vi

If a structural group $G=\{1\}$ is trivial then, up to isomorphism,
there exist exactly two Rees matrix semigroup $S=\M(\{1\};2,2;P)$
with $\Q_0 S$ a unitary ring. In the following remark, we exhibit
these semigroups since they appear as factors of a principal series
of $S$ when $\Q S$  has the hyperbolic property and contains
nilpotent elements.

\vi

\begin{obs}\label{em2}
There are four possibilities  for the Rees semigroup
$\M(\{1\};2,2;P)$, where $P$ is invertible, whose  elements are the
elementary matrices $e_{11},e_{12},e_{21},e_{22}$ and the null
matrix $\theta$ with the semigroup operation $\cdot$. Let $H=\langle
\sigma \rangle$ be the group generated by the transposition
$\sigma=(12)$ and $U=\sum_{i,j=1}^{2} e_{ij}$.
\begin{description}
\item[$\bullet$] The semigroup $M^\phi = \M(\{1\};2,2;I_{\phi}), \phi \in H $ such that if $\phi$ is the identity then $I_{\phi}$ is
the identity matrix and if $\phi =\sigma$ then $I_{\phi}=e_{12}+e_{21}$. Note
that in $M^\phi$, $e_{ij} \cdot
  e_{kl}:=e_{il}\delta_{j\phi(k)}$.
\item[$\bullet$] The semigroup $M^n=\M(\{1\};2,2;I^{n})$ with $I^n=U-e_{\sigma(n)n}$, $n=1,2$.  Note
that in $M^n$, $e_{in} \cdot e_{\sigma(n)j}=e_{ij}$, and $e_{ij}
\cdot e_{kl}:=e_{il}\delta_{jk}$ otherwise.
\end{description}

If $\phi$ is the identity of $H$ then we denote $M = M^\phi$.
Clearly the maps
$$\begin{array}{lll} \begin{array}{llll}\varphi:&M & \longrightarrow & M^\sigma\\
                             &e_{ij}& \mapsto &e_{i\sigma(j)}\end{array}& \quad \quad \quad and \quad \quad \quad&
                             \begin{array}{llll}\phi:&M^1 & \longrightarrow & M^2\\
                             &e_{ij}& \mapsto &e_{\sigma(i)\sigma(j)},  \end{array} \end{array}$$
                             are isomorphims, hence $M^\sigma \cong M {\not \cong }M^1 \cong
                             M^2$ although  $\Q_{0} M \cong \Q_{0} M^1\cong
                             M_2(\Q)$.
Clearly, the Rees matrix semigroup $M$ and $M^1$
  are generated by nilpotent elements.
\end{obs}

\vi

In the sequel we shall make free use of the following results:
\begin{enumerate}
 \item Every periodic $0$-simple semigroup (in particular any finite semigroup)   is
 completely $0$-simple, \cite[Corollary 2.56]{clfrd}. Hence, by Rees Theorem, a $0$-simple semigroups is isomorphic to some Rees matrix semigroup.
\item Let $S$ be a finite simple semigroup, if $KS$ is semi-simple then $S$ is a group, \cite[Corollary $5.24$]{clfrd}.
\item   $\Q S$ is semi-simple if, and only if,   $\Q
(S_{i}/S_{i+1})$  is semi-simple for  each factor of $S$,
\cite[Theorem $5.14$]{clfrd}.
\item  Let $\Q S$ be semi-simple. If $S_{i}/S_{i+1}$ is a factor
of $S$ then $S_{i}/S_{i+1}$ is isomorphic to a Rees  matrix
semigroup.
\end{enumerate}

\vi

 Let $S$ be a finite semigroup,
 $a,b \in S$ are inverses if $aba=a$ and $bab=b$. An inverse semigroup is a semigroup whose non-zero elements have a
unique inverse. Suppose $\Z S$ has an identity, $\U(\Z S)$  is
finite if, and only if, $S$ is an inverse semigroup which is the
disjoint union of groups that are finite Abelian groups of exponent
dividing $4$ or $6$ or $2$-Hamiltonian groups \cite[Theorem
$6.1$]{gdf}. Clearly, for such semigroups the hyperbolic property
holds.

\vi

 We shall now start a classification of all finite  semigroup whose
 semigroup algebra over $\Q$ has the hyperbolic property. In what follows we suppose that $\Z S$ has a unit.
Recall that  $S^{\theta}$ is
 nilpotent if there exists $ n \in \Z^{+}$, such that,
                   $S^{n}=\{\theta\}$. If $s \in S$ and  $s^{m}=\theta$, then
$s$ is called $m$-nilpotent element, or nilpotent. We use the expression
``{\it  nilpotent free}" to indicate the absence of non-trivial
nilpotent elements.

\vi

\begin{lemma} \label{smp}
Let $S$ be a finite semigroup. Then $\Q S$ is nilpotent free if, and
only if, $S$ admits a principal series whose factors are isomorphic
to  maximal subgroups $G$, say, of $S$ and $\Q G$ is nilpotent free.
In particular,  $S$ is the disjoint union of its maximal subgroups.
\end{lemma}
\begin{pf}
It is a consequence of \cite[Lemma $5.17$]{clfrd} and Lemma
\ref{GD}. $\square$ \end{pf}

\begin{theorem} \label{tsmp}
The algebra $\Q S$ is nilpotent free and has the hyperbolic property
if, and only if, $S$ admits a principal series for which every
factor is isomorphic to one of the groups below:
\begin{enumerate}
\item  A Higman group;
\item One of the following cyclic groups: $C_{5},C_{8}$ or $C_{12}$.
\end{enumerate}
Furthermore, at most one of the groups of type $(ii)$ occurs.
Moreover, $S$ is an inverse semigroup and it is the disjoint union
of  groups of type $(i)$ with at most one group of type $(ii)$.
\end{theorem}

\begin{pf}
Since $\Q S$ is nilpotent free, by the previous lemma, $S$ has a
principal series $S=S_1 \supset S_{2} \supset \cdots \supset
S_{n+1}=\emptyset$ whose factors $S_i / S_{i+1} \cong G_i$, a group,
and $S \cong \cup G_i$. Thus $\Q S^{\theta} \cong (\oplus_{i} \Q_0
G_{i}) \oplus \Q \theta$ and $\Gamma =(\prod \Z_{0}G_i)\times \Z
\theta$ is an order of $\Q S^{\theta}$. If $\Q S$ has the hyperbolic
property, by Theorem \ref{tfadfh} item $(i)$, $\Q S \cong \oplus
\A_i$ where at most one component $\A_{i_0}$ admits a $\Z$-order
$\Gamma_{i_0}$ such that the group $\U(\Gamma_{i_0})$ is hyperbolic
infinite. Hence, by \cite[Theorem $6.1$]{gdf}, the groups $G_i,i
\neq i_0$ are Higman groups and, by \cite[Theorem $3$]{jpp},
$G_{i_0} \in \{C_5,C_8,C_{12}\}$. Obviously, $\U(\Z S^{\theta})
\cong_{\phi} (\prod \U(\Z_{0}G_i))\times \U(\Z \theta)$, where
$\phi$ is an isomorphism. If $\theta \neq x \in S$ then $x \in G$
one of the maximal subgroup of $S$, suppose $a,b \in S$ are inverses
of $x$, then $xax=xbx=x \Rightarrow \phi(xax)=\phi(xbx)$ hence $a=b$
and $S$ is an inverse semigroup.

Conversely, if $S$ is a semigroup with a principal series whose
factors
 $S_{i}/S_{i+1} \cong G_{i}$ then $\Q_{0} S \cong \oplus
\Q_{0}( S_{i}/S_{i+1}) \cong \oplus \Q_0 G_{i}$. Consider the order
$\Gamma$ previously defined. By hypothesis, we have at most a unique
 cyclic group $G_{i_0}$, say, of order $5,8$ or $12$ and all other
 $\U(\Z G_{i}),i \neq i_0$ are trivial.  Therefore, by Theorem \ref{tfadfh} item $(i)$, $\Q S$
 has the hyperbolic property. $\square$
\end{pf}

An algebra $\A$ with the hyperbolic property and which has nilpotent
elements may be semi-simple or not. If it is semi-simple then, by
 Theorem \ref{tfadfh}, its Wedderburn-Mal\^cev decomposition has
a unique component isomorphic to $M_{2}(\Q)$. For any other
component
 the unit group of every $\Z$-order of this component is a finite group. In the next theorem we
 classify the finite semigroups whose rational semigroup algebra has these properties.

\vp

\begin{theorem}\label{tsmpn}
Let $\Q S$ be a unitary algebra with nilpotent elements. Then $\Q S$
is semi-simple and has the hyperbolic property if, and only if, $S$
has a principal series with all factors, except for one, isomorphic
to Higman groups. The exceptional one is isomorphic to a semigroup
$K$ of the following type:
\begin{enumerate}
\item $K\in \{S_{3},D_{4},Q_{12},C_{4}\rtimes C_{4}: C_4 \textrm{ acts non
    trivially on }C_4\}$;
\item  $ K\in  \{\M(\{1\};2,2;I_{d})=M \quad \quad ,\quad
  \quad\M(\{1\};2,2;
  \left(\begin{array}{ll}1&1\\0&1\end{array}\right))=M^1\}$.
\end{enumerate}
In particular, if $K$ is a group then $S$ is the disjoint union of
Higman groups and $K$.
\end{theorem}

\begin{pf}

If $\Q S$ is semi-simple then $\Q_{0}(S_{i}/S_{i+1})\cong
M_{n_{i}}(\Q G_{i})$. If $\Q S$ has the hyperbolic property then, by
 Theorem \ref{tfadfh} item $(ii)$, $\Q S \cong (\oplus \A_i)
\oplus M_2(\Q)$. Thus, for a unique index $i_0$, either $\Q_0
(S_{i_{0}}/S_{i_{0}+1}) \cong M_2(\Q)$ or $M_2(\Q)$ is the unique matrix
epimorphic image of $\Q_{0}(S_{i_0}/S_{i_0+1})$. As a consequence,
either $S_{i_0}/S_{i_0+1} \cong \m(\{1\};2,2;P)$ or
$S_{i_0}/S_{i_0+1} \cong K \in \{S_3, D_4, Q_{12}, C_4 \rtimes
C_4\}$. For each $i\neq i_0$, $n_i=1$ and $S_{i}/S_{i+1} \cong G_i$
a Higman group.

Conversely, if $S$ has a principal series as described then $\Q_0
S^\theta \cong (\oplus_{i\neq i_0} \Q_0 (S_i /S_{i+1})) \oplus
\Q_0(S_{i_0}/S_{i_0+1})$, where $G_i \cong S_i /S_{i+1}$ is a Higman
group for every $i \neq i_0$ and $K \cong S_{i_0}/S_{i_0+1}$ the
exceptional factor. Since either $\Q_0 K \cong (\oplus B_i) \oplus
M_2(\Q)$ or $\Q_0 K \cong M_2(\Q)$ we have
$$\Q_0 S \cong (\oplus \Q G_i) \oplus ((\oplus B_i) \oplus M_2(\Q)) = (\oplus
\A_i) \oplus M_2(\Q).$$ Clearly, if $\Gamma_i$ is a $\Z$-order of $\A_i$ then
$\U(\Gamma_i)$ is finite. Thus by Theorem \ref{tfadfh} item (ii)
the algebra $\Q S$ has the hyperbolic property.
$\square$
\end{pf}

\begin{proposition} \label{tri}
Let $S$ be a finite semigroup, such that, $J(\Q S) =\Q\langle
j_{0}\rangle$ for some ${j_{0}
 \in \Q S}$, and $j_{0}^{2}=0$. Then for each $s \in S$, there exists
$\lambda_{s}, \rho_{s} \in \{-1,0,1\}$ such that
 $s \cdot j_{0}=\lambda_{s}j_{0}$ and $j_{0} \cdot s=\rho_{s}j_{0}$.
\end{proposition}

\begin{pf}
  Since $J(\Q S):=J=\Q \langle j_{0}\rangle$ is an ideal, for $s \in S$
 there exists $\lambda_{s} \in \Q$, such that, $s\cdot j_{0}=\lambda_{s}
 j_{0}$. For the cyclic semigroup generated by
 $s$ there exists $n \in \N$  such that
   $e:=s^{n}$ is an idempotent \cite[$\S 1.6$)]{clfrd}. Inductively, we obtain that
 $ s^{k}\cdot j_{0}=\lambda_{s}^{k} j_{0}$ and by
 Proposition \ref{pcent}$.(ii)$ it follows that $\lambda_{s} \in
 \{-1,0,1\}$. In a similar way we obtain that $\rho_{s} \in \{-1,0,1\}$.
$\square$
\end{pf}

If  $S$ is a finite semigroup which is non-semi-simple then,
according to \cite[Corollary $5.15$]{clfrd}, every principal series
of $S$ admits a null factor (a null semigroup with two elements).

\vi
\begin{lemma} \label{nicee}
Let $S$ be a  finite non-semi-simple semigroup. If $\Q S$ has the
hyperbolic property then $S$ has a principal series with a unique
null factor of $S$. Moreover, if $S_{i_0}/S_{i_0+1}:=\{\theta,f\}$
is the null factor then either $f$ is a nilpotent of
$\hat{S}:=S/S_{i_0+1}$ and $f$ is not a nilpotent of $S$, or $f$ is
a nilpotent of $S$, $\I:=\{\theta,f\}$ is an ideal of $S$, $S/\I$ is
semi-simple and $J=\Q \langle f \rangle$. In each case, $f$ is the
unique nilpotent element of the semigroup.
\end{lemma}

\begin{pf}
Let $S^\theta =S_1 \supset S_2 \supset \cdots \supset S_n=\{\theta\}
\supset \emptyset$ be a principal series of $S$. We have for  each
factor $S_i/S_{i+1}$ and $J(\Q_0(S_i/S_{i+1}))$, the radical of
$\Q_0(S_i/S_{i+1})$, that
$$(\Q S)/J(\Q S) \cong \oplus
(\Q_0(S_i/S_{i+1})/J(\Q_0(S_i/S_{i+1}))).$$ Thus, if
$S_{i_0}/S_{i_0+1}$ is a null factor then $\Q_0 (S_{i_0}/S_{i_0+1})
\subseteq J(\Q S)$. Suppose that $S_{i_1}/S_{i_1+1}$ is another null
factor of $S$, clearly if $x_l \in S_{i_l}/S_{i_l+1}, l=0,1$, then
$\langle 1+x_0,1+x_1 \rangle \cong \Z^2$, which is contrary to
hyperbolic property of $\Q S$.

If $f$ is not nilpotent in $S$ then
  $f^2 \in S_{i_0+1} \setminus \{\theta\}$ and $f$ is nilpotent in $\hat{S}:=S/S_{i_0+1}$. If $f$ is nilpotent in $S$, by
Remark \ref{more}, $f \in J$ and consequently $J
  =\Q \langle f \rangle$. We claim that $\I:=\{\theta,f\}$ is
  an ideal of $S$. In fact, if $s \in S$ then $s f \in
  J$ and hence, by the previous proposition, $s f=\theta$ or $s f=f$ and so
  $s f \in \{\theta,f \}$. Similarly we have that $f s \in \{\theta,f
  \}$. Since in each case, $f \in S_{i_0}/S_{i_0+1}$ the unique null
  factor of the principal series of $S$, clearly $f$ is the unique nilpotent element.
$\square$
\end{pf}

\begin{proposition}\label{cst}
 Let $S$ be a finite semigroup that admits a nilpotent element $j_{0} \in S$.  $\Q S$ is non-semi-simple
 and has  the hyperbolic property if, and only if, $\I= \{\theta,
 j_{0}\}$ is an ideal of $S$ and  $S/ \I$ has a principal series whose
 factors are isomorphic to Higman groups. In particular,
  $S/ \I$ is the disjoint union of its maximal subgroups.
\end{proposition}

\begin{pf}
We have that $\Q S \cong \s(\Q S) \oplus J$ with non-trivial $J$.
Since $\Q S$ has the hyperbolic property we have, by Theorem
\ref{tfadfh}, that $\Q S \cong (\oplus \A_{i}) \oplus X$, where  ${X
\in \{J, T_{2}(\Q)\}}$   depending on the centrality of $J$. In
both cases, if $\Gamma$ is a {$\Z$-order} in $\Q S/J$, then $\U
(\Gamma)$ is finite. Therefore $\Q S/J$ has the hyperbolic property
and is nilpotent free.
By hypothesis $j_{0} \in S$ is nilpotent hence, by the previous lemma,
 ${\I :=\{j_{0}, \theta\}}$ is an
ideal of $S$ and $J= \Q \langle j_{0} \rangle$.  We have that $\Q \I
\cong \Q \langle j_{0} \rangle$ and hence ${\Q S/J \cong \Q S/ \Q \I
\cong \Q_{0}(S/ \I)}$ has the hyperbolic property and is nilpotent
free. It follows, by Theorem \ref{smp}, that $S/ \I$ admits a series
whose factors are Higman groups or the cyclic groups $C_{5},C_{8}$
and $C_{12}$. Thus the cyclic groups $C_{5},C_{8}$ and $C_{12}$ do
not occur since, by the last paragraph, $\U(\Gamma)$ is finite.

Conversely, $\I=\{j_{0}, \theta\}$ is an ideal of $S$ and $S/ \I$ admits a
series whose
 factors are Higman groups then, by   Lemma \ref{smp}, $
\Q_0 (S/ \I) \cong \oplus_{i=1}^{N}\Q G_{i}$ and  hence $\Q S/\Q \I
\cong  \oplus_{i=1}^{N}\Q G_{i}$. Since $\Q \I\cong \langle j_{0}
\rangle_{\Q}=J$, we have that the Wedderburn-Mal$\check{c}ev$
decomposition is $ \Q S \cong (\oplus \Q G_{i}) \oplus \langle j_{0}
\rangle_{\Q}$, where $\s(\Q S) \cong \oplus \Q G_{i}$ is the
semi-simple subalgebra of $\Q S$.
 If  $J$ is non-central then, by   Proposition \ref{pcent}, there exist
unique ${E_{1}, E_{N} \in \E(\Q S)}$ such that \linebreak
$E_{1}j_0=j_0 E_N=j_0$,  $j_0 E_{1}= E_Nj_0=0$ and if
$\{E_{1},E_{N},j_{0}, \theta\} \cong T_2$  then $(\Q S) (E_1+E_N)=\Q
\langle E_1,E_2,j_0 \rangle \cong \Q T_2 \cong T_2(\Q)$. Hence $\Q S
\cong \oplus_{1<i<N} \s(\Q S) E_i \oplus T_2 (\Q) \cong B \oplus
T_{2}(\Q)$, where  $B$ is in the annihilator of $\langle j_{0}
\rangle$ in $\s(\Q S)$. Furthermore, for $\A_{i}:=\s(\Q S) E_{i}$,
we have that  $ \Q S \cong( \oplus_{1<i<N} \A_{i}) \oplus
T_{2}(\Q)$. If $\Gamma \subset \Q S$ is a $\Z$-order, $\U(\Gamma)
\cong (\prod \U(\Gamma_i)) \times \U(T_2(\Z))$ where $\Gamma_i
\subset \A_i$ is a $\Z$-order. Since $ \oplus_{1<i<N} \A_i=
\oplus_{1<i<N} \s(\Q S) E_{i} \subset \oplus \Q G_{i}$ and every
$G_{i}$ is a Higman group, the unit group $\U( \Z G_{i})$ is finite
thus $\U(\Gamma_i)$ is finite. If $J$ is central then ${\Q S
\cong(\oplus \A_{i}) \oplus J}$. In both cases, the $\A_{i}$'s are
division rings and $\U(\Gamma_{i})$ is finite for any
 $\Z$-order $\Gamma_{i} \subset \A_{i}$. Thus, by Theorem \ref{tfadfh} item
 $(iii)$ and $(iv)$,
 $\Q S$ has the
 hyperbolic property.
 $\square$
 \end{pf}

\begin{theorem}\label{tnsmp}
Let $S$ be a finite semigroup. $\Q S$ is non-semi-simple
 and has  the hyperbolic property if, and only if, $S$ has a principal series with a
 unique null factor and
 the other factors are isomorphic to Higman groups.
 \end{theorem}
\begin{pf}
By Lemma \ref{nicee}, if $\Q S$ is non-semi-simple and has the
hyperbolic property then $S$ has a principal series with a unique
null factor $S_{i_0}/S_{i_0+1}=\{f,0\}:=\I$, say. If $f$ is
nilpotent in $S$ then the result follows by the last proposition.
Otherwise, since $S_{i_0+1}$ is an ideal of $S$ let
$\hat{S}:=S/S_{i_0+1}$. Then $\Q S \cong \Q S_{i_0+1} \oplus \Q_0
\hat{S}$ is a direct sum as ideals and $\Q S_{i_0+1}$ and $\Q_0
\hat{S}$ has the hyperbolic property. Clearly, $\hat{S}$ has the
nilpotent $f$ and by the last proposition the factors of $\hat{S}$,
and thus the factors $S_i/S_{i+1}, 1 \leq i <i_0$, are isomorphic to
Higman groups. If $\Gamma$ is an order of $\Q_0 \hat{S}$ then
$\U(\Gamma)$ is virtually cyclic hence $\Q S_{i_0+1}$ is nilpotent
free and has the hyperbolic property. Thus by  Theorem \ref{tsmp}
the factors of $S_{i_0+1}$, and therefore the factors $S_i/S_{i+1}, i_0
<i$, are Higman groups.

Conversely, on the conditions over the factors of a series of $S$ we
have that \linebreak  $\Q S \cong \Q S_{i_0+1} \oplus \Q_0 \hat{S}$
is a direct sum as ideals. By Theorem \ref{tsmp}, $\Q S_{i_0+1}$ has
the hyperbolic property and is nilpotent free. By the last
proposition $\Q_0 \hat{S}$ is hyperbolic. Clearly, If $\Gamma$ is an
order of $\Q S$ then $\U(\Gamma)$ is hyperbolic and the result now
follows. $\square$
\end{pf}

 In \cite{jpsf}, \cite{js} and \cite{thss} are classified the quadratic extensions
 $K=\Q [\sqrt{-d}]$, where $-d$ is a square free integer, and the
finite groups $G$ for which the  group ring $\mathfrak o_K[G]$ of
$G$ over the ring of integers of $K$ has the
property that the group $\mathcal U_1(\mathfrak o_K[G])$ of units of
augmentation 1 is hyperbolic. Therefore it is
 natural to classify the extensions $K$ and the semigroups $S$, such that the
 algebra $K S$ has the hyperbolic property. By remark \ref{hphi}, if $KS$
 has nilpotent elements then, since the integral basis of  $o_K$ has two
 elements, $KS$ does not have the hyperbolic property. Therefore a
 necessary condition for $KS$ to have the hyperbolic property is that $KS$ must
 be nilpotent free.

\vp

\begin{theorem}
Let $K=\Q(\sqrt{-d})$, where $-1 \neq d \neq 0$, and $S$ a finite
semigroup. The algebra $KS$ is nilpotent free and has the hyperbolic
property if and only if $KS$ is semi-simple, $S$ admits a principal
series whose factors are isomorphic to Higman groups $G$ or at most one copy of
a group $H$, where $G, H$ and $d$ are determined by the following table.
$$\begin{tabular}{|c|c|c|}
\hline
Higman groups $G$ & Group $H$&$K=\Q(\sqrt{-d})$\\
\hline & & $0<d$\\
\cline{2-3}       &$C_3$     &$0<d \neq 3$\\
\cline{2-3}       an elementary Abelian $2$-group      & $C_4$     &$0<d\neq 1$\\
\cline{2-3}        &$Q_8$     &$d \equiv 7 \pmod 8,0<d$\\
\hline an Abelian group of
exponent dividing $4$      &$C_8$        & $d=1$\\
\hline an Abelian group of exponent dividing $6$ & & $d=3$\\
\hline
\end{tabular}$$
Moreover, $S$ is an inverse semigroup and it is a disjoint union
of  groups $G$ with, if it is the case, a unique group $H$.
\end{theorem}

\begin{pf}
If $K S$ is nilpotent free and has the hyperbolic property then so
is $\Q S$ and, by Theorem \ref{tsmp}, $S$ has a principal series whose
factors are Higman groups and at most one of the groups $C_5, C_8,
C_{12}$. By \cite[Theorem $4.7$]{jpsf}, only $C_8$ is such that
$\U(\mathfrak o_K C_8)$ is hyperbolic and in this case $d=1$. If a factor
of $S$ is an Abelian Higman group $G$ then the free rank of
$\U(\mathfrak o_K G)$ is $0$ or $1$, since $\U(\Gamma)$ is
hyperbolic for all $\Z$-order $\Gamma \subset KS$. Using \cite[Theorem
$3.7$]{jpsf} we determine the groups $G$ and $H$ and the extension
$K$. If $G$ is non-Abelian then, by \cite[Theorem $4.7$]{jpsf},
$G=Q_8$, $0<d$ and $d \equiv 7 \pmod 8$. Since $\U(\mathfrak o_K
Q_8)$ is an infinite hyperbolic group, the Higman group allowed as a
factor is only an elementary Abelian $2$-group, see \cite[Corollary
$3.5$]{jpsf}. Conversely, if $\Gamma \subset KS$ is a $\Z$-order it
follows immediately by \cite{jpsf}  that $\U(\Gamma)$ is hyperbolic.
Clearly, $S$ is an inverse semigroup.
\end{pf}

\section{Idempotents of Maximal Subgroups}

\vi

Let $S$ be a finite semigroup with a nilpotent element $j_0$. In
this section, we investigate the idempotents of maximal subgroups of
$S$. In case $\Q S$ is non-semi-simple and has the hyperbolic
property, the study of idempotents enable us to obtain more
information on the structure of $S$. In fact
we prove in the last theorem that $S$ has some explicit semigroups as {\it basic
blocks} of its structure which we define below as $T_2,\hat{T}_2 $ and $T'_2$.

\vi

\begin{definition}
As a set $T_2=\hat{T}_2=\{e_1,e_2,j_0,\theta\}$ and $T'_2=T_2 \cup
\{e_3\}$ are semigroups with the operation $\cdot$ given by the Cayley
table:

$$\begin{array}{lll}
T_2&\hat{T}_2&T'_2\\
\begin{tabular}{r|r|r|r|r}
$\cdot$ &$e_1$&$e_2$&$j_0$&$\theta$\\
\hline
$e_1$&$e_1$&$\theta$&$j_0$&$\theta$\\
\hline
$e_2$&$\theta$&$e_2$&$\theta$&$\theta$\\
\hline
$j_0$&$\theta$&$j_0$&$\theta$&$\theta$\\
\hline $\theta$&$\theta$&$\theta$&$\theta$&$\theta$
\end{tabular} \qquad \qquad &
\begin{tabular}{r|r|r|r|r}
$\cdot$ &$e_1$&$e_2$&$j_0$&$\theta$\\
\hline
$e_1$&$e_1$&$j_0$&$j_0$&$\theta$\\
\hline
$e_2$&$\theta$&$e_2$&$\theta$&$\theta$\\
\hline
$j_0$&$\theta$&$j_0$&$\theta$&$\theta$\\
\hline $\theta$&$\theta$&$\theta$&$\theta$&$\theta$
\end{tabular} \qquad \qquad&
\begin{tabular}{r|r|r|r|r|r}
$\cdot$& $e_1$&$ e_2$&$ e_3$& $j_0$& $\theta$\\
\hline
$e_1$& $e_1$ &$e_3$ & $e_3$ &$j_0$& $\theta$\\
\hline
$e_2$&$e_3$ &$e_2$ & $e_3$& $\theta$& $\theta$\\
\hline
$e_3$& $e_3$& $e_3$ & $e_3$ & $\theta$ &$\theta$\\
\hline
$j_0$&$\theta$ &$j_0$&$\theta$& $\theta$& $\theta$\\
\hline $\theta$&$\theta$& $\theta$& $\theta$& $\theta$& $\theta$
\end{tabular}
\end{array}$$
\end{definition}
\vi

In what follows $S=\cup G_{i}\cup \{\theta,j_{0}\},
j_0^2=\theta$, see Proposition \ref{cst}, and $N=|\E(\Q S)|>2$.
If $E_l \in \E(\Q S)$ then $ 1=\sum_{1<l<N} E_{l} +E_1+E_N$. Let
$E:=E_{1}+E_{N}$ and  $e \in \Q S$ be any idempotent, hence $e=
\sum_{1<l<N} eE_{l}+eE$, where $(eE_{l})^{2}=eE_{l}\in \A_l$ which
is, by Theorem \ref{tfadfh}, a division ring, $\forall \ 1<l<N$. Therefore,  $eE_{l} \in \{E_{l},0\}$. Let $E_{e_l}:=eE_l \neq 0$ thus $e= \sum
E_{e_{l}} + eE$.

\vi

\begin{proposition}
If $e_{i}$ is the group identity element of the group $G_i$ then
$e_{i}$ has one of the following expressions: $$\begin{array}{l}\sum
E_{i_{l}}+ E_{1}+\lambda j_{0}\\\sum
 E_{i_{l}}+  E_{N}+\mu j_{0}\\ \sum E_{i_{l}}+  E \\ \sum
 E_{i_{l}} \end{array}$$ with $0 \neq E_{i_l}:=e_i E_l, E_l \in \E(\Q S)$ and for some $\lambda, \mu \in \Q$,
 Moreover  the last two expressions are central idempotents.
\end{proposition}

\begin{pf}

Write   $e_{i}=\sum E_{i_{l}}+ u E_{1}+v E_{N}+w j_{0}$ (recall that
the $E_{i_{l}}$ are orthogonal, central, they \linebreak annihilate
$j_{0}$, $E_{1}j_0=j_0 E_N=j_0$ and $E_{N}j_0=j_0 E_1=0$).  Hence
\linebreak $e_{i}^{2}=\sum E_{i_{l}}+ u^{2} E_{1}+v^{2} E_{N}+w(u+v)
j_0=e_i$ and thus $u,v \in \{1,0\}$ and $w(u+v)=w$. If $u=v=1$ then
$w=0$ and therefore $e_{i}=\sum E_{i_{l}}+E_{1}+E_{N}$. The others
possibilities are : $u=1,v=0,w=1$, $u=0,v=1,w=1$, and $u=v=w=0$,
resulting in the other expressions. $\square$ \end{pf}

\vi

\begin{lemma}\label{gidem}
Suppose $\Q S\cong (\oplus \A_{i}) \oplus \Q j_0$, $E_i \in \E(\oplus \A_i)$ and
 $E_{1}j_{0}=j_{0}E_{N}=j_{0}$. Let $G,H$ be maximal subgroups of $\ S$. If  $\A_{1}
 \subseteq \Q G$ then  $gj_{0}=j_{0},\ j_{0}
 g=0, \ \forall g\in G$ and if $\A_{N} \subseteq \Q H$ then
 $j_{0}h=j_{0}, \ h j_{0}=0, \forall h\in H$.
\end{lemma}

\begin{pf}
We have that $E_1 \in \A_1 \subset \Q G$, and so $E_1=\sum_{g \in G}
\alpha_g g$. By the property of $E_1$ it holds $0\neq j_0=E_1
j_0=(\sum \alpha_{g}\lambda_g) j_0$ and, by Proposition \ref{tri},  the  $\lambda_g 's \in
\{0, \pm 1\}, \forall g\in G$. Therefore, there exists $g_0 \in G$ such
that $\lambda_{g_0}=1$. If $e_1$ is the identity of $G$ then it follows that $ e_1 j_0=j_0$ and thus
$g j_0=j_0, \forall g\in G$, because $G$ is a finite group and
$\{\theta  , j_0\}$ is an ideal. Similarly, if $e_N$ is the identity of $H$ then
$j_0 e_N=j_0$ and $j_0 h=h$ for all $h\in H$.

 Since $\Q j_0$ is an ideal we have that $j_0 e_1=\rho j_{0}$ and $\rho \in
 \{0, \pm 1\}$. Suppose $0 \neq \rho =1$, say, thus $e_1 j_0=j_0=j_0 e_1$; then  $e_1$
 centralizes $j_0$ and hence $e_1 \notin \A_{1}$, a contradiction.
 In the same way we prove that $e_N j_0=0$. $\square$
 \end{pf}

\begin{proposition}
Let $G$ be a maximal subgroup of $S$. Denote by $e \in G$ its
identity element and suppose that $e j_0=j_0$. If $e=\sum E_{e_l}+
E_{1}+\lambda j_{0}$ then, $\forall g \in G$, \linebreak $g=\sum
gE_{e_l}+ E_{1}+\lambda j_{0}$. Also  if  $e=\sum E_{e_l}+
E_{N}+\mu j_{0}$ then, $\forall g \in G$, $g=\sum g E_{e_l}+
E_{N}+\mu j_{0}$.
\end{proposition}

\begin{pf}
For  $g \in G$, we have that $g=ge=\sum gE_{e_l}+ g E_{1}+g\lambda
j_{0}$.  By   Lemma  \ref{gidem}, $g j_{0}=j_0$ and hence ${g=\sum
gE_{e_l}+ g E_{1}+\lambda j_{0}}$.

To determine $g E_1$, recall that $\Q \langle E_1,E_N,j_0\rangle$ is an ideal
 of $\Q S$. So we may
 write \linebreak  $g E_1 =t E_1 +s E_N + r j_0$. There exists $k \in \Z$ such that $g^{k}=e$; since the orthogonality of
 $E_1$ respect $E_l,l\neq 1$ and $E_{1}j_{0}=j_{0}$ we conclude that
 $E_1 g^{k}=E_1+\lambda j_0$. By comparing with the equation $(g E_1)^{k} =(t E_1 +s
 E_N + r j_0)^{k}=t^{k}E_1+s^{k}E_{n}+r'j_{0}$ we reach: $t^{k}=0, s^{k}=1$ and
 $r'=\lambda$. Hence, $g=\sum g E_{e_l}\pm E_{1}+\lambda j_{0}$ and, multiplying on the right by
 $j_0$ and using Lemma \ref{gidem}, we obtain that $g=\sum g E_{e_l}+ E_{1}+\lambda j_{0}$.

For the other case: $e=\sum E_{e_l}+  E_{N}+\mu j_{0}$, it holds that
$j_0 e=j_0$. If $g \in G$ then similarly $g=\sum g E_{e_l}+
E_{N}+\mu j_{0}$. $\square$
\end{pf}

\begin{theorem}\label{tesf}
 Let  $e_{1} \in G_{1}$ and $e_{N} \in G_{N}$ be the group identities and suppose that $e_1 j_0 =j_0 e_N=j_0$. Write
 $$\begin{array}{l}e_1=\sum E_{1_{l}}+  E_{1}+\lambda j_{0}\\
                  e_N=\sum E_{N_{l}}+  E_{N}+\mu j_{0}. \end{array}$$
Then only one of following options holds:
\begin{description}
\item[$(i)$] $$e_l e_k =0 \Leftrightarrow e_{k}e_{l}=|\alpha|j_0 \textrm{ and } \lambda+\mu=\alpha, \textrm{ where } \{1,N\}= \{l,k\}$$
and either $\alpha=0$ and $\{ e_{1},e_{N},j_0,\theta\} \cong T_2$ or
$\alpha \in \{-1,1\}$ and $\{ e_{1},e_{N},j_0,\theta\} \cong
\hat{T}_{2}$.
\item[$(iii)$] $$\textrm{If } e_N e_1 \neq 0 \textrm{ then } e_1 e_N = e_N e_1=:e_3
  \textrm{ and } \lambda+\mu=0,$$
and $\{e_1,e_N,e_3,j_0,\theta\} \cong T'_2 $.
\end{description}
\end{theorem}

\begin{pf}
Since the idempotents $E_{i} \in \E(\Q S)$ are orthogonal,
$j_0e_1=j_0E_1=e_n j_0=E_N j_0=0$, and $e_1j_0=E_1j_0=j_0e_n
=j_0E_N=j_0$
\begin{eqnarray}
 \nonumber  e_1 e_N&=&\sum E_{1{_l}}E_{N_{l}}+(\lambda +\mu)j_{0}\\
 \nonumber  e_N e_1&=&\sum E_{N{_l}}E_{1_{l}} \\
   \label{qc}         e_1e_N&=&e_Ne_1 +(\lambda +\mu)j_{0}
\end{eqnarray}

Without loss of generality suppose $l=1,k=N$.

If $e_1 e_N=0$ then $-(\lambda+\mu)j_0=e_N e_1$. In addition, if
$\lambda+\mu=0$ then $e_N e_1=0$, and the converse is clear. If
$\lambda+\mu \neq 0$ then $e_N e_1$ is a non-trivial nilpotent
element of $S$. Thus $e_N e_1=j_0$, since $S$ has a unique nilpotent
element, clearly $\lambda+\mu=-1=\alpha$ (for $l=N,k=1$ we have
$\alpha=1$). The converse is straightforward.

If $e_1 e_N$ and $e_N e_1$ are non-zero elements then, by equation
(\ref{qc}), the set $\{e_1 e_N, e_N e_1,j_0\} \subseteq S$ is
$\Q-$L.D. Since any element of this set is not zero thus $e_1
e_N=e_N e_1:=e_3$ is a non-trivial idempotent and $\lambda+\mu=0$.
The converse is clear.
$\square$
\end{pf}

The semigroups $T_2$, $T'_2$ e $\hat{T}_2$ are, in some sense, the
basic  building  blocks of the semigroups $S$  whose rational
semigroup algebra is non-semi-simple and has the hyperbolic
property.

\vp

\section*{Acknowledgements}
This article corresponds to the second chapter of the third author's
PhD. Thesis supervised by the second author to whom he is grateful. We are indebted to the referee whose suggestions considerably improved the manuscript.

\vp

\end{document}